\documentclass[10pt]{amsart}
\usepackage{latexsym, amsmath,amssymb}
 \numberwithin{equation}{section}
 \newcommand{\half}{{\textstyle \frac 12}}

\newcommand{\ccal}{\mathcal{C}}

\newcommand{\lcal}{\mathcal{L}}

\newcommand{\Q}{{\mathbb Q}}
\newtheorem{theo}{{\sc Theorem}}[section]

\newtheorem{lem}[theo]{{\sc Lemma}}

\newenvironment{rem}{\medskip\noindent{\it Remark:\/} }{\medskip}

\newtheorem{defn}[theo]{{\sc Definition}}

\renewcommand{\epsilon}{\varepsilon}
\newtheorem{theorem}{Theorem}

\newtheorem{lemma}[theorem]{Lemma}
\newtheorem{corr}[theorem]{Corollary}

\newtheorem{proposition}[theorem]{Proposition}
\newtheorem{deff}[theorem]{Definition}

\newcommand{\bth}{\begin{theorem}}
\newcommand{\ble}{\begin{lemma}}
\newcommand{\bcor}{\begin{corr}}

\newcommand{\bdeff}{\begin{deff}}

\newcommand{\bprop}{\begin{proposition}}
\newcommand{\ele}{\end{lemma}}
\newcommand{\ecor}{\end{corr}}
\newcommand{\edeff}{\end{deff}}

\numberwithin{theorem}{section}

\newcommand{\eprop}{\end{proposition}}

\newcommand{\Rn}{{\mathbb R}^n}

\newcommand{\la}{\lambda}

\newcommand{\e}{\varepsilon}

\newcommand{\supp}{\text{supp }}
\renewcommand{\Pi}{\varPi}

\renewcommand{\epsilon}{\varepsilon}

\newcommand{\R}{{\mathbb R}}
\newcommand{\Z}{{\mathbb Z}}

\newcommand{\1}{{\rm 1\hspace*{-0.4ex}%
\rule{0.1ex}{1.52ex}\hspace*{0.2ex}}}

\begin{document}

\title[Focal points and sup norms of eigenfunctions]
{Focal points and sup-norms of eigenfunctions 
}

\author[C. D. Sogge]{Christopher D. Sogge}
\address{Department of Mathematics, Johns Hopkins University, Baltimore, MD 21218, USA}
\email{sogge@jhu.edu}
\author{Steve Zelditch}
\address{Department of Mathematics, Northwestern University, Evanston, IL 60208, USA}
\email{s-zelditch@northwestern.edu}

\thanks{Th research of the first author was partially supported by NSF grant \# DMS-1069175 and the second by  \#  DMS-1206527.  The first author was also supported in part
by a Simons Fellowship.}

\begin{abstract} If $(M,g)$ is a compact real analytic Riemannian manifold, we give a necessary and sufficient condition for there to be a sequence of quasimodes of order $o(\la)$  saturating sup-norm estimates. In particular,
it gives optimal conditions for existence of eigenfunctions satisfying maximal sup norm bounds.  The condition is that there exists a self-focal point $x_0\in M$ for the geodesic flow at which the associated Perron-Frobenius operator
$U_{x_0}: L^2(S_{x_0}^*M) \to L^2(S_{x_0}^*M)$ has a nontrivial invariant $L^2$ function.  The proof is based on an explict 
Duistermaat-Guillemin-Safarov pre-trace formula and  von Neumann's ergodic theorem.
\end{abstract}

\maketitle

\section{Introduction and main results}

Let $(M,g)$ be a compact boundaryless Riemannian manifold.  We then let 
$\Phi_t(x,\xi)=(x(t),\xi(t))$ denote the  homogeneous Hamilton flow on $T^* M \backslash 0$  associated to the principal symbol $|\xi|_g$
of $\sqrt{-\Delta}$, with $\Delta=\Delta_g$ denoting the Laplace-Beltrami operator of $(M, g)$. Since  $\Phi_t$ 
preserves the unit cosphere bundle $S^*M = \{|\xi|_g = 1\}$ , it defines a flow on  $S^*M$ which preserves  the Liouville measure.  For a given $x\in M$, let ${\mathcal L}_x\subset S^*_xM$ denote those unit directions $\xi$ for which $\Phi_t(x,\xi)\in S_x^*M$ for some time $t\ne 0$, and let 
$|{\mathcal L}_x|$ denote its surface  measure $d\mu_x$ in $S^*_x M$  induced by the Euclidean metric $g_x$.  Thus, ${\mathcal L}_x$ denotes the initial directions of geodesic loops through $x$.

In \cite{SZ} it was shown that
\begin{equation}\label{1.1}
\|e_\la\|_{L^\infty(M)}=o(\la^{\frac{n-1}2}),
\end{equation}
if 
$$|{\mathcal L}_x|= 0 \quad \text{for all } \, x\in M,$$
whenever $e_\la$ is an $L^2$-normalized eigenfunction of frequency $\la$, i.e.
\begin{equation}\label{1.2}
(\Delta+\la^2)e_\la = 0, \quad \text{and } \, \, \int_M|e_\la|^2 \, dV=1.
\end{equation}
Here, $dV=dV_g$ denotes the volume element associated with the metric.
Moreover, there is a corresponding bound for the $L^2\to L^\infty$ norm of projection operators 
 onto shrinking spectral bounds, i.e.,
\begin{equation}\label{1.o}
\bigl\|\chi_{[\la,\la+o(1)]}\bigr\|_{L^2(M)\to L^\infty(M)}=o(\la^{\frac{n-1}2}).
\end{equation}
By this we mean that, given $\e>0$, we can find a $\delta(\e)>0$ and $\Lambda_\e<\infty$ so that
\begin{equation}\label{1.o2}
\bigl\|\chi_{[\la,\la+\delta(\e)]}f\|_{L^\infty(M)}\le \e \la^{\frac{n-1}2} \, \|f\|_{L^2(M)}, 
\quad \la \ge \Lambda_\e,
\end{equation}
with $\chi_{[\la,\la+\delta]}:L^2(M)\to L^2(M)$ denoting projection onto frequencies
(i.e., eigenvalues of $\sqrt{-\Delta}$) in the interval $[\la, \la+\delta]$.  

In \cite{STZ} this result was improved by showing that if ${\mathcal R}_x\subset{\mathcal L}_x\subset S^*_xM$ is the 
set of recurrent directions over $x$ and if $|{\mathcal R}_x|=0$ for all $x\in M$, then
we have \eqref{1.1}.


We shall now give a simple and natural further refinement in the real analytic case which involves an ergodicity condition.  The real analytic case is simple to analyze due to the fact that if $x\in M$
then there are just two extreme possibilities regarding the nature of 
${\mathcal L}_x\subset S^*_xM$, the loop directions.  Indeed, as shown in \cite{SZ}, 
either $|{\mathcal L}_x|=0$ or ${\mathcal L}_x=S^*_xM$.  In the second case there is also a minimal time $\ell>0$ so that $\Phi_\ell(x,\xi)\in S^*_xM$, meaning that all geodesics starting at $x$ loop back at exactly this minimal time $\ell$. We call
such a point a {\em self-focal point}. 
We note that there may exist more than one self-focal point and the   minimal common return time $\ell$ of the loops may depend on $x$ but for simplicity of notation
we do not give it a subscript.  If we then write
\begin{equation}\label{1.5}
\Phi_\ell(x,\xi)=(x,\eta_x(\xi)), \quad \xi\in S^*_xM,
\end{equation}
then the {\em first return map}, 
\begin{equation} \label{FRM} \eta_x: S^*_xM\to S^*_xM \end{equation} above our self-focal point is real analytic.  Following Safarov~\cite{Sa}, we can associate to this first return map
the {\em  Perron-Frobenius operator}
$U_x:L^2(S^*_xM, d\mu_x)\to L^2(S^*_xM, d\mu_x)$
by setting
\begin{equation}\label{1.6}
U_xf(\xi)=f(\eta_x(\xi))\, \sqrt{J_x(\xi)}, \quad f\in L^2(S^*_xM, d\mu_x),
\end{equation}
where  $J_x(\xi)$ denotes the Jacobian of the first return map, i.e. $\eta_x^* d\mu_x = J_x(\xi) d\mu_x$.   Clearly $U_x$
is a unitary operator and \begin{equation} \label{eta} \eta_x^* (f d\mu_x) = U_x(f) d\mu_x.  \end{equation}

The key assumption underlying our results is contained in the following:
\bigskip

\begin{defn}\label{DIS} 
A self-focal point $x\in M$ is said to be {\em dissipative} if $U_x$ has no invariant  function $f \in L^2(S^*_x M)$.
Equivalently, $\Phi_{\ell}$ has no invariant $L^1$ measure with respect to $d\mu_x$.
\end{defn}
The dissipative condition is a spectral condition on $U_x$. If $U_x $ has any $L^2$ eigenfunction $g$ then
$U_x g = e^{i \theta} g$, and, since $U_x$ is a positive operator, $U_x |g| = |U_x g| = |g|$. Hence the dissipative
condition is the condition that the spectrum of $U_x$ is purely continuous. For this reason, one might prefer the
term `weak mixing' ; but that  might create the wrong impression that $\Phi_{\ell}$ is weak mixing with respect to some given invariant measure.
The term `dissipative'  refers to the Hopf decomposition of $\Phi_{\ell}$ on $S^*_x M$ into conservative and dissipative parts.
As discussed in   \S \ref{HOPF}, lack of an $L^2$ eigenfunction does not necessarily imply that the conservative part is of
measure zero, so the term `dissipative point' does not precisely mean that $\Phi_x$ is a dissipative dynamical system.


Our main result then is the following

\begin{theorem}\label{theorem1.1}  Let $(M,g)$ be a real analytic compact boundaryless manifold of dimension $n\ge2$.  Then 
\eqref{1.1} and \eqref{1.o}  are true if and only if every self-focal point is dissipative. 
\end{theorem}

\begin{rem}  The condition that every self-focal point is dissipative  is also equivalent to 
\begin{equation}\label{1.1b}
\lambda^{-1} \|\nabla e_\la\|_{L^\infty(M)}=o(\la^{\frac{n-1}2}).
\end{equation}
For expository simplicity we will not prove it here, but the proof is almost the same.
\end{rem}

On any surface of revolution,  the
invariant eigenfuntions have sup norms of order $\simeq \la^{\frac{n-1}2} $ and achieve their suprema at the poles. Their gradients
achieve their suprema at a distance $\frac{C}{\lambda}$ away from the poles. The above condition rules out such
examples. 

Theorem \ref{theorem1.1} is stronger than the  result in \cite{STZ}.  Indeed,  if $|{\mathcal R}_x| = 0$
then  $U_x$ cannot have an invariant $L^2$ function and so \eqref{1.o} would hold.  The reason is that  
existence of $f \in L^2$ with  $U_x f = f$ implies  (by \eqref{eta}) that   $|f|^2 d\mu_x$
is  a finite invariant measure in the class of $d\mu_x$.   By the Poincar\'e recurrence theorem,  $|f|^2d\mu_x$ -
almost every point would be recurrent, and thus the set of recurrent points would have positive $\mu_x$ measure.

Theorem \ref{theorem1.1} is in a sense the optimal result on the problem of relating sup-norms of eigenfunctions
to properties of the geodesic flow. As will be discussed at the end of the introduction,  the result can only be sharpened
by further understanding of the   dynamical question of relating existence of an invariant $f \in L^2(S^*_x M, d\mu_x)$ to the structure
of loops of $(M, g)$. 

\subsection{Geometry of  loops and  self-focal points}

Associated to a  self-focal point $x$  is the flowout  manifold
\begin{equation} \label{LAMBDA} \Lambda_{x} = \bigcup_{0 \leq t \leq \ell} G^t S^*_{x} M.  \end{equation}
It is an immersed Lagrangian submanifold of $S^*M$ whose projection
$$\pi : \Lambda_{x} \to M $$
has a ``blow-down singuarity" at $t = 0, t = \ell$. For this reason self-focal points were called blow-down points in \cite{STZ}.
We may view $\Lambda_{x_0}$ as the embedding of  
the mapping cylinder  ${\mathcal C}_x$ of $\eta_x$, i.e. as 
\begin{equation} {\mathcal C}_x = S_x^*M \times [0, \ell] /
\cong,\;\;\;\mbox{where}\;\; (\xi, \ell) \cong (\Phi_{\ell}(x, \xi), 0)
\end{equation} In \cite{STZ}  it is proved that  the map
$$\iota_x(\xi, t) = G^t (x, \xi): {\mathcal C}_x \to \Lambda_x \subset  S^*_x M$$
 is a Lagrange immersion whose image is
$\Lambda_{x} \subset S^*_x M$. If $\ell$ is the minimal period of all loops (i.e. if
there are no exceptionally short loops) then $\iota_z|_{S^*_z
\times (0, \ell)}$ is an embedding.

Focal points come in two basic kinds, depending on  the
first return map $\eta_x$.   We say that $x$  is  a {\it pole}  if  
 $$ \eta_x = Id_x: S^*_x M \to S^*_x M. $$   Equivalently, the set $\ccal \lcal_x$ of smoothly closed geodesics
based at $x$ is all of $S^*_x M$, 
\begin{equation} \label{POLE}  x \;\mbox{is a pole if} \;\;\ccal \lcal_x = S^*_x M, \;\;\;\; \ccal \lcal_x = \{\xi \in S^*_x M: \Phi_t(x, \xi) = (x, \xi)\}.
\end{equation}
On the other hand, it is possible that $\eta_x = Id$ only on  a set of zero measure in $\lcal_x$, which in the 
analytic case means that it is almost nowhere the identity and $\lcal_x$ must have  codimension $\geq 1$.  We call such a $\eta_x$ {\it twisted}.  Thus,
\begin{equation} \label{TWISTED} x  \; \mbox{is self-focal with a twisted return map if} \;\; \mbox{codim} \;(\ccal \lcal_x) \geq 1. 
\end{equation}

  Examples of poles are 
the poles $x$ of a surface of revolution (in which case all geodesic loops at $x$ are smoothly closed). Examples
of self-focal points with twisted return map are the 
four umbilic points of two-dimensional tri-axial ellipsoids, from  which  all geodesics loop back 
 at time $2 \pi$ but are almost never smoothly closed \cite{K}. The only smoothly closed directions 
are the geodesic (and its time reversal) defined by the middle length `equator'. 
There are topological restrictions on manifolds possessing a self-focal point. 
In \cite{BB} a manifold
with such a point is denoted a $F_{\ell}^{x_0}$ (or $Y^{x_0}_{\ell}$-)-manifold; if $\ell$ is the least common return time for all loops
it is denoted by $L^{x_0}_{\ell}$.  
If  $(M, g)$ has a focal point $x_0$ from which  all geodesics are simple
(non-intersecting) loops, then the integral cohomology ring $H^*(M, \Z)$  is generated by one element  \cite{N}. 
For  an $F_{\ell}^{x_0}$ manifold, $H^*(M, \Q)$ has a single generator (Theorem 4 of \cite{BB}).  Most results
on manifolds with self-focal points consider only the special case of Zoll metrics; see \cite{Be} for
classic results and \cite{Ol}  for some recent results
and references.

In the case of a triaxial ellipsoid  $E \subset \R^3$, the first return map $\eta_x$ is a totally dissipative  expanding map of the circle
with  two   fixed points, one 
a source and one a sink. It has invariant $\delta$-measures at the fixed points and an  infinite  locally
$L^1$ invariant measure on each component of the  complement. According to Theorem \ref{theorem1.1}, the
eigenfunctions of $E$ cannot achieve maximial sup norm bounds. In fact, the result of \cite{STZ} already rules
out maximial eigenfunction growth on the ellipsoid.

Note also that the analog of \eqref{1.o} is automatic in the Euclidean case since,
by Plancherel's theorem and duality,
$$\Bigl\| \, \int_{\{\xi\in \Rn: \, |\xi|\in [\la,\la+\delta]}
e^{ix\cdot \xi} \Hat f(\xi)\, d\xi \, \Bigr\|_{L^\infty(\Rn)}
\le c_n\sqrt{\delta} \la^{\frac{n-1}2}\|f\|_{L^2(\Rn)}, \quad \la\ge 1.$$
Thus, by Theorem~\ref{theorem1.1}, the condition as to whether or not there
is a self-focal point $x\in M$ so that $U_xg=g$ for some $0\ne g\in L^2(S^*_xM)$, with
$U_x$ as in \eqref{1.6}, provides a necessary and sufficient condition determining
when the spectral projection operators on $M$ also enjoy improved sup-norm bounds
over shrinking intervals.  Besides \cite{SZ} and \cite{STZ}, earlier works on
related problems are in \cite{BS}, \cite{BSSY}, \cite{Toh}, \cite{China} and 
\cite{Stein}.  We should also point out that in \eqref{1.o2} one needs that
$\delta(\e)\to 0$ as $\e\to 0$ since, for any compact $n$-dimensional Riemannian
manifold $(M,g)$, by (5.1.12) in \cite{Sbook}, one has 
$$\limsup_{\la\to \infty}\la^{-\frac{n-1}2}\|\chi_{[\la,\la+1]}\|_{L^2(M)\to L^\infty(M)}>0.$$

\subsection{Coherent states associated to self-focal points}

The proof of Theorem \ref{theorem1.1} extends beyond eigenfunctions to certain
kinds of quasi-modes, which play a fundamental role in the proof. 
To construct them, we fix $\rho\in {\mathcal S}(\R)$ satisfying
\begin{equation}\label{1.8} 
\rho\ge 0, \, \, \,  \rho(0)=1, \quad 0\le \hat \rho\le 1, \, \, \text{and } \, \,
\Hat \rho(t)=0 \, \, \text{if } \, t\notin (-1,1).
\end{equation}
We set 
$$P=\sqrt{-\Delta},$$
and consider the  operators 
\begin{equation} \label{rhoa} \rho(T(\la-P)) = \frac{1}{T} \int_{\R} \hat{\rho}(\frac{t}{T})  e^{i t \lambda} e^{- it P} dt. \end{equation}

If we freeze the second component of the Schwartz kernel of \eqref{rhoa} at a focal point $x_0$ we obtain 
a  semi-classical Lagrangian quasi-mode 
\begin{equation} \label{qm} \psi^{x_0}_{\lambda, T}(x) :=  \rho (T(\lambda - \sqrt{\Delta})(x, x_0)  \end{equation}
associated to the Lagrangian submanifold \eqref{LAMBDA} in the sense of \cite{CdV, D}. That is,
it is a semi-classical oscillatory integral with large parameter $\lambda$ whose phase generates \eqref{LAMBDA}.
More precisely, $\{\psi^{x_0}_{\lambda, T}\}$ is a one-parameter family of quasi-modes depending on the parameter
$T$ (as well as the semi-classical parameter $\lambda$).  We  refer to \eqref{qm}  
 as  {\it coherent states} centered at $x_0$.  In terms of   an  orthonormal basis of
real valued eigenfunctions   $\{e_j\}_{j=0}^\infty$, 
$$\psi^{x_0}_{\lambda, T} (x) = \la^{-\frac{n-1}2} \, \sum_{j=0}^\infty \rho\bigl(T(\la-\la_j)\bigr) e_j(x) e_j(x_0), $$

The following Lemma explains the precise sense in which \eqref{qm} are quasi-modes:

\begin{lem} For each $T$, $\psi^{x_0}_{\lambda, T}$ is a semi-classical Lagrangian distribution associated
to the Lagrangian \eqref{LAMBDA} whose  principal symbol pulls back under $\iota_x$ to  $$\frac{1}{T} e^{i \lambda t} \hat{\rho}(\frac{t}{T}) |dt|^{\half} \otimes |d\mu_x|^{\half}.$$
One has
\begin{equation}\label{v.16}
\|\psi_{\la,T}^{x_0}\|_{L^2}\le C,
\end{equation}
and, moreover, for any $\epsilon$ there exists $T_0$ so that  for $T \geq T_0$,
\begin{equation}\label{v.17}
||(\Delta + \lambda^2) \psi_{\lambda, T}^{x_0} ||_{L^2}  \leq C \epsilon \lambda. \end{equation}
 \end{lem}
For fixed $T$, \eqref{qm} is therefore a  very kind of weak quasi-mode
but as a family it behaves like  a quasi-mode of order $o(\lambda)$ as $T \to \infty$.

Further, the  value of $\psi_{\lambda, T}^{x_0}$ at its peak point $x_0$ equals
\begin{equation}\label{1.9}
\la^{-\frac{n-1}2} \, \sum_{j=0}^\infty \rho\bigl(T(\la-\la_j)\bigr) (e_j(x_0))^2,
\end{equation}
and as we will see, this peak value is of maximal growth if and only if $U_x$ \eqref{1.6} has an invariant $L^2$
function.

The fact that $\psi_{\lambda,T}$ is a semi-classical Lagrangian quasi-mode is an immediate consequence 
of the well-known  parametrix construction for $e^{-i t P}$ (see \eqref{parametrix}), from which it is
clear  that one may express \eqref{qm} as an oscillatory
integral  with large parameter $\lambda$. The symbol can also be evaluated as in the Lemma directly from this
expression and the principal symbol of $e^{-it P}$ computed in \cite{DG}. Note that the symbol is invariant
under the geodesic flow on the long time interval where $\hat{\rho}(\frac{t}{T}) \equiv 1$. As $T \to \infty$
one gets a kind of time average over orbits and we will see that the family \eqref{qm} behaves as a family
of approximate quasi-modes with invariant symbols. It is natural to ask whether the invariant $L^2$ function
can be `quantized' to construct a quasi-mode of order $o(\lambda)$; we discuss this in \S \ref{FUTURE}.

It is also not hard to see that \eqref{v.17} is valid.  Indeed, since $\rho\in {\mathcal S}(\R)$, we have for any $N=1,2,\dots$,
\begin{align*}
\int_M \bigl|(\Delta+\la^2)\psi^{x_0}_{\la,T}(x)\bigr|^2 \, dV&=\la^{-(n-1)}\sum_j (\la^2-\la_j^2)^2 \, \bigl(\rho(T(\la-\la_j))\bigr)^2 (e_j(x_0))^2
\\
&\le C_N\la^{-(n-1)}\sum_j (\la^2-\la_j^2)^2 \, (1+T|\la-\la_j|)^{-N} (e_j(x_0))^2
\\
&\le C_NT^{-2}\la^{-(n-1)}\sum_j (\la+\la_j)^2\, (1+T|\la-\la_j|)^{-N} (e_j(x_0))^2
\\
&\le C_N  T^{-2} \, \la^2,
\end{align*}
if $N-2>n$, using in the last step the fact that
$$\sum_{\la_j\in [\mu,\mu+1]}(e_j(x))^2\le C(1+\mu)^{n-1}, \quad \mu\ge 0.$$
Using this fact one also easily obtains \eqref{v.16}.

\subsection{Outline of the proof of Theore \ref{theorem1.1}} 

We now outline  the strategy we shall employ in
proving Theorem \ref{theorem1.1} and explain the role of the quasi-modes \eqref{qm}.

First, we shall  prove the $o(\la^{\frac{n-1}2})$ bounds \eqref{1.1} and \eqref{1.o} under the assumption that
every self-focal point is dissipative. 
As we pointed out before, because of our real analyticity assumption, there are two extreme cases: {\em self-focal points} and {\em non-focal points} (where $|{\mathcal L}_x|=0$).
Thus, it is natural to split our estimates for eigenfunctions or quasi-modes into two cases.  The first, which is new, is to prove
favorable bounds in a ball around  each dissipative self-focal point. The second, which was dealt with in \cite{SZ}, is to prove such bounds near a non-focal point.

After we prove \eqref{1.o}, we shall prove the converse direction:
the $o(\la^{\frac{n-1}2})$ bounds at self-focal points for which $U_x$ \eqref{1.6}  has no nonzero invariant $L^2$ functions can be turned around 
to show that there are $\Omega(\la^{\frac{n-1}2})$ bounds if  nontrivial
invariant functions do exist. 

The main estimate is as follows: If  all focal points are dissipative, then  for any $\e>0$, we prove that there are finite
$T_0=T_0(\e)$ and  $\Lambda_0=\Lambda_0(\e)$
so that (as in \eqref{1.9})
\begin{equation}\label{1.10}\sum_{j=0}^\infty \rho\bigl(T_0(\la-\la_j)\bigr)(e_j(y))^2\le
\e\la^{n-1} \quad \la\ge \Lambda_0, \, \, y\in M,
\end{equation}
Since $\rho(0)=1$, it is plain that this estimate implies the first assertion,
\eqref{1.1}, of the theorem.
As  mentioned above, we split the proof of \eqref{1.10}  into two very different cases: (i)    \eqref{1.10} is valid in a neighborhood of a given dissipative self-focal point,  and (ii)
\eqref{1.10} is valid in a neighborhood of any non-focal point.


The role of the quasi-modes \eqref{qm} in the proof of the first part of Theorem~\ref{theorem1.1}  is explained by 
the following

\begin{proposition}\label{prop1.2}  Let $x$ be a dissipative self-focal point in our
real analytic compact boundaryless manifold of dimension $n\ge 2$.  Then for any $\e>0$,
we can find a neighborhood ${\mathcal N}(x,\e)$ of $x$ and finite numbers 
$T=T(x,\e)$
and  $\Lambda=\Lambda(x,\e)$ so that
\begin{equation}\label{1.11}
\sum_{j=0}^\infty
\rho\bigl(T(\la-\la_j)\bigr)(e_j(y))^2 \le \e \la^{n-1}, \quad
\text{if } \, 
 y\in {\mathcal N}(x,\e), \, \,
\text{and } \, \, \la\ge \Lambda.
\end{equation}
\end{proposition}

The proof of the estimate \eqref{1.11} at a self-focal point $y$ uses the dissipative assumption and the von Neumann mean
ergodic theorem. To extend the estimate to points $y$ in a neighborhood of a self-focal point, we use the rather
explicit formula for the wave trace formula for large times T of Safarov \cite{Sa,SV}, which shows 
the leading part of the left side of \eqref{1.11}
varies smoothly.  Consequently, if the inequality in \eqref{1.11}  holds at the self-focal point, it holds in some neighborhood of the
self-focal point when $\la$ is sufficiently large. 

\begin{rem} 
The proof of Theorem \ref{theorem1.1} and Proposition \ref{prop1.2}  shows that the suprema of $|e_\la(x)|$, resp.  $|\nabla e_\la(x)|$,  is obtained in a ball of radius $O(\frac{1}{\lambda})$ around a self-focal point. But it does not show that the
suprema are obtained at a self-focal point, and
obviously  both suprema cannot be attained  at the same point.  Since $e_\la$ oscillates on the scale $\frac{1}{\lambda}$,
it is apriori possible that it vanishes at the self-focal point and takes its suprema in a ball of radius $\frac{C}{\lambda}$
around the self-focal point. Thus,  Proposition \ref{prop1.2}  shows that 
not only is
$|e_{\la}(y)|$ small at the dissipative self-focal point but also that there is ``propagation of smallness'' to a neighborhood of such a point.
\end{rem}

The other case, which we need to handle in order to prove the first part of Theorem~\ref{theorem1.1}, is to show that we also have these types of bounds near every non-focal point:

\begin{proposition}\label{prop1.3}  Let $x\in M$ and assume that $|{\mathcal L}_x|=0$.
Then, given any $\e>0$, we can find a neighborhood ${\mathcal N}(x,\e)$ of $x$ and finite numbers $T=T(x,\e)$ and $\Lambda=\Lambda(x,\e)$ so that
\eqref{1.11} is valid.
\end{proposition}

 Proposition~\ref{prop1.3} is a special case of Theorem~1.1 in \cite{SZ}, but we shall present its simple proof for the sake of completeness.  As in \cite{SZ}, it uses an idea coming from Ivrii's~\cite{Iv} proof of his generalization of the Duistermaat-Guillemin theorem \cite{DG}.

Before outlining the proofs of the Propositions, let us see how they imply
\eqref{1.o}.  
We first claim that \eqref{1.10} 
is valid if and only if,
given $\e>0$, there are numbers $\delta_0=\delta_0(\e)$ and $\Lambda_0=\Lambda_0(\e)$ so that for every $y\in M$
\begin{equation}\label{1.12}
\sum_{\{j: \, |\la-\la_j|\le \delta_0\}} \, \bigl(e_j(y)\bigr)^2
\le \e \la^{n-1}, \quad \text{if } \, \,
\la\ge \Lambda_0.
\end{equation}
This clearly implies \eqref{1.o} since
$$\bigl\|\chi_{[\la,\la+\delta]}\bigr\|_{L^2(M)\to L^\infty(M)}^2
=\sup_{y\in M}\sum_{\{j: \, \la_j\in [\la, \la+\delta]} \bigr(e_j(y)\bigr)^2.$$
We shall also make use of this fact when we prove the first part of Corollary~\ref{corr}. 


To prove our claim
that \eqref{1.10} is equivalent to \eqref{1.12},
we first note that, by \eqref{1.8}, $\rho(s)\ge 1/2$ when
$|s|\le \delta$, for some fixed $\delta>0$.  From this it is clear that \eqref{1.10} implies \eqref{1.12}.
Since $\rho\in {\mathcal S}(\R)$, it
follows that for any $N=1,2,3,\dots$ there is a constant $C_N$ so that
$$\sum_{j=0}^\infty \rho(T_0(\la-\la_j)) (e_j(y))^2
\le C_N \sum_{j=0}^\infty \bigl(1+T_0|\la-\la_j|\bigr)^{-N}(e_j(y))^2.$$
As a result if we take $\delta_0$ in \eqref{1.12} to be $1/T_0$ and
$N=n+1$, we conclude that \eqref{1.12} implies \eqref{1.10}.

We now use a  compactness argument to show   that the  Propositions imply \eqref{1.12}. 
 Propositions~\ref{prop1.2} 
 implies that,
if $\e>0$ is given and if $x$ is a self-focal point 
then there must be a neighborhood ${\mathcal N}(x,\e)$ of $x$ and numbers
$\delta_x>0$ and $\Lambda_x<\infty$ so that
\begin{equation}\label{1.13}
\sum_{|\la_j-\la|\le \delta_x}\, \bigl(e_j(y)\bigr)^2
\le \e \la^{n-1}, \quad
\text{if } \, y\in {\mathcal N}(x,\e), \, \, \text{and } \, 
\la \ge \Lambda_x.
\end{equation}
 Proposition~\ref{prop1.3} implies the same conclusion
for any non-focal point $x$.  
Since $\{{\mathcal N}(x,\e)\}_{x\in M}$ is an open covering of $M$, by Heine-Borel, there must be a finite subcovering.   In other words, there must be points
$x_j\in M$, $1\le j\le N$, so that
$M\subset \bigcup_{j=1}^N{\mathcal N}(x_j,\e)$.
By \eqref{1.13}, if 
$\Lambda_0=\max\{ \Lambda_{x_1}, \dots, \Lambda_{x_N}\}$,
and 
$\delta_0=\min\{\delta_{x_1},\dots, \delta_{x_N}\},
$
then we must have \eqref{1.12}. As discussed above, this implies \eqref{1.10} , hence \eqref{1.o} and therefore
\eqref{1.1}.

\subsection{Generalization to fat quasi-modes}

With a small  additional effort, the proof of Theorem \ref{theorem1.1}  establishes 
 a more general  result for growth rates of  ``fat quasi-modes" or quasi-modes of order $o(\lambda)$. They are defined as follows:
When  $n=2$ and $n=3$, we say that  
a sequence $\{\phi_{\la_k}\}$  of {\em quasi-modes of order $o(\lambda)$} if
\begin{equation}\label{q.1}\
\int_M|\phi_{\la_k}|^2 \, dV=1, \quad \text{and } \, \|(\Delta+\la^2_k)\phi_{\la_k}\|_{L^2(M)}=o(\la_k).
\end{equation}
In higher dimensions, $n\ge 4$, as explained in \cite{STZ}, pp. 164-165, one needs
to modify this definition in order to get natural results by requiring that 
\begin{equation}\label{q.2}
\int |\phi_{\la_k}|^2 \, dV=1, \quad
\text{and } \, \, \|S_{[2\la_k,\infty)}\phi_\la \|_{L^\infty(M)}+
\|(\Delta+\la^2)\phi_{\la_k}\|_{L^2(M)}=o(\la_k),
\end{equation}
if $S_{[2\la_k,\infty)}:L^2(M)\to L^2(M)$ denotes the projection onto the $[2\la_k,\infty)$ frequencies.  As was shown in \cite{STZ}, condition \eqref{q.1} automatically implies \eqref{q.2} if
$n=2$ or $n=3$, and, moreover, if $|{\mathcal R}_x|=0$ for each $x\in M$, then
\begin{equation}\label{q.3}
\|\phi_\la \|_{L^\infty(M)}=o(\la^{\frac{n-1}2}),
\end{equation}
whenever $\phi_\la$ is a sequence of quasi-modes of order $o(\lambda)$.

As we shall see, an equivalent formulation of Theorem~\ref{theorem1.1} is the following.

\begin{corr}\label{corr}
Let $(M,g)$ be a compact real analytic boundaryless manifold of dimension $n\ge 2$.
Then \eqref{q.3} holds for quasi-modes of order $o(\lambda)$ if and only if 
every self-focal point of $M$ is dissipative.
\end{corr}

We shall prove this easy consequence of Theorem~\ref{theorem1.1} in the final section of the paper and also state a natural problem about whether the above quasi-mode
condition can be weakened if the sup-norms are $\Omega(\la^{\frac{n-1}2})$.


\subsection{\label{HOPF} Hopf decomposition and existence of a finite invariant measure in the class of $d\mu_x$}
In this section we explain the term `dissipative' in Definition \ref{DIS}. 

Let $\Phi$ be an invertible measurable map of a measure space $(X, \mu)$. 
A set $W$ is called wandering if the sets $\Phi^{-k} W$ for $k
\geq 0$ are disjoint, i.e. if no point of $W$ returns to $W$.
$\Phi$ is called {\it conservative}  if there exists no wandering set of
positive measure.

 $\Phi$ is
called {\it recurrent } if for all Borel sets $A$, almost every point of
$A$ belongs to the set $A_{rec}$ of points returning at least once to $A$,
$$A_{rec} := A \cap \bigcup_{k = 1}^{\infty} \Phi^{-k} (A). $$
It is called infinitely recurrent if almost all points of $A$
belong to
$$A_{inf} := \{x \in A: \Phi^k \in A, \;\; \mbox{for infinitely many }\; k \geq 1\} = A \cap \bigcap_{n = 1}^{\infty}
\bigcup_{k = n}^{\infty} \Phi^{-k} (A). $$
 The recurrence theorem (see \cite{K}, Theorem
3.1) states that the following are equivalent: \begin{itemize}

\item $\Phi$ is conservative;

\item $\Phi$ is recurrent;

\item $\Phi$ is infinitely recurrent;


\end{itemize}

The Hopf decomposition states:

\begin{theo} (cf. \cite{K}, Theorem 3.2)  If $\Phi$ is null-preserving and non-singular, then
there exists a decomposition of $X$ into two
disjoint measurable sets $C, D$ (the conservative and dissipative
parts), so that

\begin{itemize}

\item $C$ is $\Phi$-invariant;

\item $\Phi |_C$ is conservative;

\item There exists a wandering set $W$ so that $D =
\bigcup_{-\infty}^{\infty} \Phi^k W. $

\end{itemize}

\end{theo}

 We are of course interested in the Hopf decompsition for  $\Phi = \Phi_{\ell},
X = S^*_x M$ and $\mu = \mu_x$. 
The term `dissipative' in Definition \ref{DIS} is intended to suggest that $(\Phi_{\ell}, S^*_x M, d\mu_x)$
is dissipative if there exists no invariant $ f \in L^2(S^*_x M, d\mu_x) $  for the unitary operator $U_x$. 
It is obvious that if such an invariant $f$ exists, then $\mu_x(C) > 0$.  However, in the class of measurable
dynamical systems, there  exist examples (due to Halmos) of  conservative systems without finite invariant measures.
The recent review \cite{EHIP} contains  several equivalent criteria for existence of a finite invariant measure in
the class of $\mu$. 
It is not clear to us whether $\mu_x(C) > 0$ in our setting implies existence of an invariant $L^2$ function, and
we do not claim that $\Phi_{\ell}$ is dissipative if $x$ is a dissipative self-focal point.

\subsection{\label{FUTURE} Potential further improvements}

There are three potential avenues for improvement of Theorem \ref{theorem1.1}. 

The first is purely geometric or dynamical. It is possible that existence of an invariant $L^2$ function
for $U_x$ at some self-focal point $x$ implies that $x$ is a pole.  Furthermore, it is possible that
there do not even exist real analytic Riemannian manifolds $(M, g)$ with $\dim M \geq 3$ which have self-focal
points with twisted first return maps. 
 For instance, it  does not appear that ellipsoids of
dimension $\geq 3$ have self-focal points. Having asked a number of experts on loops, we have the impression
that this question is open.  

It would also simplify the proof if we knew that there are only a finite number of self-focal points with
twisted first return maps for analytic $(M, g)$. It is not hard to prove that if there does exist an infinite number of
such points, then the first return times tend to infinity as they approach a limit point.  All known examples
have only a finite number of self-focal points with twisted first return maps. For further comments on
the geometry we refer to \cite{Z}. 

Second, it is a natural question if one can use an invariant $L^2$ function $f$ in $L^2(S^*_x M)$ to construct
a better quasi-mode. The classical construction of a  quasi-mode satisfying
$$||(\Delta + \lambda^2) \psi_{\lambda} ||_{L^2}  \leq C$$
as an oscillatory integral \cite{D} requires a smooth invariant density on \eqref{LAMBDA}. On the one hand,
it is possible that the invariant $f$ is always $C^{\infty}$. On the other hand, it is possible that one can
quantize a `rough' density to construct a quasi-mode. In the real analytic case, there is a possible construction
using the FBI transform which we plan to investigate in future work. 

Finally, it is natural to generalize the argument of this article to general $C^{\infty} $ metrics. Our assumption
that $g$ is real analytic is only to simplify the geometric analysis of loopsets. In general, they can be badly
behaved, and that requires further approximations and arguments. We plan to give the general proof in
a subsequent article. However, the present one contains the main analytical ideas.

\subsection{Organization}

The paper is organized as follows.  In the next section we shall prove that
we have the $o$-bounds \eqref{1.o} provided that the Perron-Frobenius operators
\eqref{1.6} above every self-focal point have no nonzero invariant functions.  Then in \S 3, we shall see that this proof can be easily modified to give the remaining
part of Theorem~\ref{theorem1.1}, by showing  that we cannot have \eqref{1.12}
if $y=x$ is a self-focal point where the operator $U$ in \eqref{1.6} has a nontrivial
$L^2(S_x^*M)$ function.  Finally, in \S 4 we shall give the simple argument showing
that Theorem~\ref{theorem1.1} implies the Corollary and discuss further problems
for quasi-modes.

\section{Proof of $o(\la^{\frac{n-1}2})$ bounds}

As we noted above, to prove the $o$-bounds posited in Theorem~\ref{theorem1.1}, it
is natural to split the analysis into two cases: bounds near self-focal points and
near non-focal points. Let us start by giving the argument for the former.

\subsection{Analysis near  self-focal points with no invariant $L^2$ functions}

In this subsection we shall prove Proposition~\ref{prop1.2}.  So let us assume that $x$ is a self-focal point, and that $\ell>0$ is the first return time for the geodesic
flow through $x$.  Thus, $\Phi_\ell(S^*_xM)=S^*_xM$, with $\ell$ being the minimal such time.  We also are assuming that the associated Perron-Frobenius $U$ has
no nontrivial invariant $L^2(S^*_xM)$ functions.

Our assumption about $\ell$ does not rule out the existence of subfocal times $0<t_0<\ell$ and directions $\xi\in S^*_xM$ with $\Phi_{t_0}(x,\xi)\in S^*_xM$ (i.e., loops through
$x$ of length shorter than $\ell$).  On the other hand, the set ${\mathcal E}_x\subset S^*_xM$ of such directions must be closed and of measure zero.  As a result, it will
be a simple matter to modify the argument for the case of no subfocal times to handle the general case.  We shall do so at the end of this  subsection.

To simplify the notation, we note that, after possibly rescaling the metric, we may, and shall take $\ell$ to be equal to one.  We then let $\eta=\eta_x: S^*_xM\to S^*_xM$, as in
\eqref{1.5}, be the first return map and $U=U_x$, as in \eqref{1.6}, the associated unitary operator on $L^2(S^*_xM)$.  Assuming that the
operator $U$ in \eqref{1.6} has no nonzero invariant functions, we must show that there is a 
neighborhood ${\mathcal N}(x)$ of $x$ in $M$ so that we have \eqref{1.11}.

Using the Fourier transform, we can rewrite the left hand side of \eqref{1.11} as
\begin{equation}\label{rho}
\sum_j \rho(T(\la-\la_j))e_j(x) \, e_j(y)=\frac1{2\pi T}\int \Hat \rho(t/T) \, \bigl(e^{itP}\bigr)(x,y) \, e^{-it\la} \, dt,\end{equation}
and, because of our assumptions, the integrand vanishes if $|t|\notin (-T,T)$.  Here
$\bigl(e^{itP}\bigr)(x,y)$
denotes the kernel of the half-wave operator $e^{itP}$.  It follows from H\"ormander's theorem~\cite{Hsing} on the propagation of
singularities and our temporary assumption that there are no subfocal times, that $(t,y)\to \bigl(e^{itP}\bigr)(y,y)$ is 
smooth when $y=x$ and $t\in \R\backslash {\mathbb Z}$.

On the other hand, $\bigl(e^{itP}\bigr)(x,x)$ will be singular when $t=\nu\in \Z$.  We recall that for $t$ near a given such
$\nu$ and $y$ near $x$, one can write down 
a parametrix for this half-wave operator which is a finite sum of terms of the form
\begin{equation}\label{2.5}
(2\pi)^{-n}\int_{\Rn} e^{iS_{\nu,j}(t,y,\xi)-iz\cdot \xi}a_{\nu,j}(t,y,z,\xi)\, d\xi, \quad j=1,\dots, N(\nu),
\end{equation}
where, up to Maslov factors, $i^{s_\nu}$,  (see \cite[p. 68]{DG}) the principal symbols of the $a_{\nu,j}$ are nonnegative zero order symbols and the phase functions $S_{\nu,j}(t,y,\xi)$ are real-valued, homogeneous of degree one in $\xi$ and generating functions
for  portions of the canonical relation associated with $P$.  That means that, if $p(x,\xi)$ is the principal symbol of $P$, then
\begin{equation}\label{b.6}
\partial_t S_{\nu,j}(t,y,\xi)=p(y,\nabla_x S_{\nu,j}(t,y,\xi))  .
\end{equation}
We therefore can write for $|t|\le T$ and $y$ in a small neighborhood of $x$
\begin{equation}\label{parametrix}
\bigl(e^{itP}\bigr)(y,y)=(2\pi)^{-n}\sum_{|\nu|\le T}\sum_{j=1}^{N(\nu)}
\int_{\Rn}e^{iS_{\nu,j}(t,y,\xi)-iy\cdot \xi}a_{\nu,j}(t,y,y,\xi)\, d\xi +O_T(1),
\end{equation}
where, because of our assumption that there are no subfocal points, we may assume that
\begin{equation}\label{parametrix2}
a_{\nu,j}(t,x,x,\xi)=0, \quad \text{if } \, \, t\notin (\nu-\delta, \nu+\delta),
\end{equation}
with $\delta>0$ small but fixed.

By the Hamilton-Jacobi equations associated with $p(x,\xi)$, this  means that if $t=\nu\ell =\nu$ and we have that,  $a_{\nu,j}(\nu, x,x,\xi)\ne 0$, for
some $\nu=0, \pm 1, \dots, \pm T$, then
\begin{equation}\label{b.7}
\nabla_x S_{\nu,j}(\nu,x,\xi)=\eta^\nu(\xi), 
\quad \xi \in S^*_{x}M,
\end{equation}
where $\eta(\xi)$ is as in \eqref{1.5}, and 
$\eta^\nu =\eta \circ \eta \circ \cdots \circ \eta$ is the $\nu$-fold composition of $\eta$ if $\nu>0$,  while $\eta^\nu=\eta^{-\nu}$ if $\eta<0$, and $\eta^0(\xi)=\xi$.  Additionally, since $e^{i\nu P} \circ (e^{i\nu P})^*=I$, one sees, by a theorem of H\"ormander (see  e.g., \cite[Theorem 6.1.4]{Sbook}) that  we must also have that for $\nu$ as above 
\begin{multline}\label{b.8}
\sum_{j=1}^{N(\nu)}|a_{\nu,j}\bigl(\nu,x,x,\eta^\nu(\xi)\bigr)|=\sqrt{J^\nu(\xi)} \quad 
\text{mod } \,S^{-1}_{1,0},
\\ J^\nu(\xi)=|\det (d\eta^\nu(\xi)/d\xi)|, \, \, \, \xi\in S^*_{x}M.
\end{multline}
Based on this, we conclude that in order to prove \eqref{1.11}, it suffices to show that, given $\e>0$, there is a neighborhood ${\mathcal N}(x,\e)$ of $x$
so that if $T$ is large and fixed, we have
\begin{multline}\label{b.9}
T^{-1}\left| \,  \sum_{|\nu|\le T}\sum_{j=1}^{N(\nu)} \la^n \iint \Hat \rho(t/T)\,  e^{i\la( S_{\nu,j}(t,y,\xi)-y\cdot \xi-t)}  \, 
a_{\nu,j}(t,y,y,\la \xi)\, d\xi dt\, \right| \le \e \la^{n-1} 
\\
+O_{\e,T}(\la^{n-2}), \quad y\in {\mathcal N}(x,\e).
\end{multline}

To prove this, as in \cite{DG} and \cite{Sa}-\cite{SV}, we shall use stationary phase in the $r$ and $t$ variables, if we write $\xi=r\omega$,
where $\omega\in S^*_yM$.  Given our assumption that
$\eta^\nu(S^*_{x}M)=S^*_{x}M$, $\nu =0, \pm 1, \pm 2, \dots$, this is all we can do.  In view of this, it is natural to rewrite \eqref{b.9} in the equivalent form
\begin{equation}\label{b.10}
T^{-1}\Bigl| \,  \sum_{|\nu|\le T}\sum_{j=1}^{N(\nu)} A_{\nu,j}(y) \Bigr| \le \e \la^{-1} 
+O_{\e,T}(\la^{-2}), \quad y\in {\mathcal N}(x,\e),
\end{equation}
where
\begin{equation}\label{b.11}
A_{\nu,j}(y)= (2\pi)^{-1}\iint \Hat \rho(t/T)\,  e^{i\la( S_{\nu,j}(t,y,\xi)-x\cdot \xi-t)}  \, a_k(t,y,y,\la \xi)\, d\xi dt.
\end{equation}

We shall get the gain $\la^{-1}$ in the right side of \eqref{b.10} from stationary phase in these two variables,
while an additional gain of $O(\e)$ at $y=x$ will come from von Neumann's ergodic theorem if $T=T(\e)$ is large and fixed.  Uniformity over a small neighborhood of $x$ (depending on $\e$ and $T$) will come from stationary phase with parameters, due to the fact that the  Hessian of the phase function in the $t, r$ variables is nondegenerate
for every fixed direction $\omega\in S^*_{x}M$
at $x$. 
We shall use the fact that if the leading coefficient in the stationary phase expansion is small at a self-focal point, then it is small in a neighborhood of this 
point.

Let 
\begin{equation}\label{b.12}
\psi_{\nu,j}(t,y,\xi)=S_{\nu,j}(t,y,\xi)-y\cdot \xi -t
\end{equation}
denote the phase function in  \eqref{b.11}.  
%
Let us first argue that \eqref{b.10} is valid when $y$ is our self-focal point $x$.  We follow an argument in \cite{Sa}, \cite{SV}.

If we fix a direction $\omega\in S^*_{x}M$ (so that $p(x,\omega)=1)$ and write
$\xi=r\omega$, then by the Euler homogeneity relations we have
\begin{equation}\label{b.13}
\partial_r\psi_{\nu,j}(t,x,r\omega)=0 \, \iff \, S_{\nu,j}(t,x,\xi)-x\cdot \xi =0.
\end{equation}
In other words
\begin{equation}\label{b.14}
\partial_r\psi_{\nu,j}(t,x,\xi)= 0 \, \implies \psi_{\nu,j}(t,x,\xi)=-t.
\end{equation}
Additionally, by \eqref{b.6},
\begin{equation}\label{b.15}
\partial_t\psi_{\nu,j}(t,y,\xi)=0\implies p(y,\nabla_x S_{\nu,j}(t,y,\xi))=1,
\end{equation}
meaning that if $a_{\nu,j}(\nu,x,x,\xi)\ne 0$ and if $\eta^\nu(\xi)=\nabla_x S_{\nu,j}(t,x_0,\xi))$ is as in \eqref{b.7}, then
$\eta^\nu(\xi)\in S^*_{x}M$.   The Hessian of $\psi_{\nu,j}$ in the $t,r$ variables therefore must be of the form
\begin{equation}\label{b.16}
H(x)=\left( \begin{array}{cc} 0 & p(x,\eta^\nu(\xi)) \\  p(x,\eta^\nu(\xi)) & 0\end{array} \right),
\end{equation}
in view of \eqref{b.6} and the fact that $p(x,\eta^\nu(\xi))$ is homogeneous
of degree one.
Thus,
\begin{equation}\label{b.17}
\det H(x)=-1.
\end{equation}
Therefore, by stationary phase (e.g., \cite[Theorem 7.7.5]{H1}), we have
\begin{multline}\label{b.18}
A_{\nu,j}(x)=(2\pi)^{-1}\int\int_0^\infty \int_{-T}^T e^{i\la\psi_{\nu,j}(t,x,r\omega)}
\hat \rho(t/T)a_{\nu,j}(t,x,x,\la r\omega) r^{n-1}dr dt d\omega
\\
= \la^{-1} \int  e^{-i\la \nu}\hat \rho(\nu/T)a_{\nu,j}(\nu,x,x,\la\omega) \, d\omega  +O_T(\la^{-2}).
\end{multline}
By \eqref{b.8},
\begin{equation}\label{b.20}
\sum_{j=1}^{N(\nu)} |a_{\nu,j}(\nu,x,x,\la\omega)|=  \sqrt{J^\nu(\omega)}
+O_T(\la^{-1}).
\end{equation}

By \eqref{b.18},  \eqref{b.20} and the fact that $0\le \hat \rho\le 1$, we have
\begin{multline}\label{b.21}
T^{-1}\sum_{|\nu|\le T}\sum_{j=1}^{N(\nu)}|A_{\nu,j}(x)|
\le \la^{-1} T^{-1}\sum_{|\nu|\le T} \sum_{j=1}^{N(\nu)}
\int_{S^*_xM}|a_{\nu,j}(\nu,x,x,\omega)| \, d\omega +O_T(\la^{-2})
\\
=\la^{-1} \,  T^{-1}\sum_{|\nu|\le T} 
\int_{S_{x}^*M}\sqrt{J^\nu(\omega)} \, d\omega  
 + \, O_T(\la^{-2}).
\end{multline}

Therefore, we would have \eqref{b.10} if we could show that, given $\e>0$, we can chose $T\gg 1$ so that
\begin{equation}\label{b.22}
T^{-1}\sum_{|\nu|\le T} \int_{S^*_{x_0}M}\sqrt{J^\nu(\omega)} \, d\omega <\e/2.
\end{equation}
If $\1$ denotes the function on $S^*_{x}M$ which is identically one, and if
$U:L^2(S^*_{x}M)\to L^2(S^*_{x}M)$ is as in \eqref{1.6}, then we can rewrite
the left hand side of \eqref{b.22} as
$$\bigl\langle \, T^{-1}\sum_{|\nu|\le T} U^\nu \1, \, \1 \, \rangle.$$
Since $U:L^2(S^*_xM)\to L^2(S^*_xM)$ is unitary, by von Neumann's ergodic theorem 
$$\frac1{2T}\sum_{|\nu|\le T} U^\nu \1 \to \Pi(\1), \quad \text{in }\, L^2(S^*_xM), \quad
\text{as } \, T\to +\infty,$$ 
where $\Pi(\1)$ denotes the projection of $\1$
onto the $U$-invariant subspace of $L^2(S^*_xM)$.  Our assumption 
that this operator has no nonzero invariant functions means that $\langle \Pi(\1), \1\rangle =0$, and therefore
we have \eqref{b.22} if $T=T(\e)$ is sufficiently large.

We have shown that 
\begin{equation}\label{b.23}
T^{-1} \sum_{|\nu|\le T}\sum_{j=1}^{N(\nu)}|A_{\nu,j}(x)| \le \e \la^{-1}/2 \, + \, O_T(\la^{-2}).
\end{equation}
Let us now argue that because of \eqref{b.17}, this implies that we can find
a neighborhood ${\mathcal{N}(x,\e)}$ of $x$ so that \eqref{b.10} is valid 
which
would finish the proof of \eqref{b.9}.

The situation at $x$ was simple since the phase functions $\psi_{\nu,j}$, as a function of $r$ and
$t$, only had critical points on the support of the amplitudes exactly at $r=1$ and $t=\nu$.  
If we assume that $y$ is sufficiently close to
$x$, and, as we may, that the  $a_{\nu,j}$ satisfy \eqref{parametrix2} with $\delta>0$ sufficiently small, it follows from \eqref{b.17} and the implicit function theorem that there 
are unique pairs $(t_{\nu,j}(y,\omega), r_{\nu,j}(y,\omega))$, 
which depend smoothly on $y$ and $\omega\in S^*_{y}M$ so that if we replace $x$
by $y$ in the oscillatory integral in \eqref{b.18}, the unique
stationary point of the phase $(t,r)\to \psi_{\nu,j}(t,y,r\omega)$ occurs at this point.  Also, assuming that $y$ is close to $x$ we will have that $r_{\nu,j}(y,\omega)$ is very close to $1$ and $t_{\nu,j}(y,\omega)$ is very close to $\nu$ if $a_{\nu,j}(\nu,x,r\omega)\ne0$.    
Let $H_{\nu,j}(y,\omega)=H_{\nu,j}(t_{\nu,j}(y,\omega),r_{\nu,j}(y,\omega))$ denote the Hessian of the phase function
at the stationary point.  It then follows that $\det H_{\nu,j}$ is close to $-1$ if $y$
is close to $x$.
Therefore by stationary phase (e.g.
\cite[Theorem 7.7.6]{H1}, for $y$ close to $x$, we have the following analog
of \eqref{b.18}
\begin{multline}\label{b.24}
A_{\nu,j}(y)=
\int \int_0^\infty \int_{-T}^T e^{i\la\psi_{\nu,j}(t,y,r\omega)}
\hat \rho(t/T)a_{\nu,j}(t,y,y,\la r\omega) r^{n-1}dr dt d\omega
\\
=\la^{-1}   
\int |\det H_{\nu,j}|^{-1/2}e^{-i\la t_{\nu,j}(y,\omega)}\hat \rho(t_{\nu,j}(y,\omega)/T)
\\ 
\qquad \qquad \qquad\times a_{\nu,j}(t_{\nu,j}(y,\omega),y,y,\la r_{\nu,j}(y,\omega)\omega)\, r_{\nu,j}(y,\omega)^{n-1} d\omega
\\ +O_T(\la^{-2}),
\end{multline}
where the constant in the $O_T(\la^{-2})$ error term can be chosen to be uniform for $y$
sufficiently close to $x$.  We conclude that
\begin{multline*}\left|\sum_{|\nu|\le T}\sum_{j=1}^{N(\nu)}\bigl(|A_{\nu,j}(y)|-|A_{\nu,j}(x)|\bigr)\right|
\\
\le \la^{-1} 
\sum_{|\nu|\le T}\sum_{j=1}^{N(\nu)}\int
 \Bigl||\det H_{\nu,j}(y,\omega)|^{-1/2} \hat  \rho(t_{\nu,j}(y,\omega)) 
 \\
\qquad  \qquad \qquad\qquad\qquad\qquad\qquad\qquad\times |a_{\nu,j}(t_{\nu,j}(y,\omega),y,y,\la r_{\nu,j}(y,\omega)\omega)| \,  r_{\nu,j}(y,\omega)^{n-1}
\\
-\hat \rho(\nu)|a_{\nu,j}(\nu, x,x,\la \omega)|\Bigr| \, d\omega 
\, + \, O_T(\la^{-2}).
\end{multline*}
Since $a_{\nu,j}\in S^0_{1,0}$ and $t_{\nu,j}(y,\omega)$ and $r_{\nu,j}(y,\omega)$ are smooth functions for $y$ near $x$ which satisfy $t_{\nu,j}(x,\omega)=\nu$ and
$r_{\nu,j}(x,\omega)=1$, it follows that each of the preceding integrals is a smooth function of $y$ near $x$, which vanishes when $y=x$.
Therefore, 
by \eqref{b.23}, we deduce that \eqref{b.10} must be
valid if ${\mathcal N}(x,\e)$ is a sufficiently small neighborhood of
$x$, since at this point, $T$ has been fixed.

\subsection*{Handling the contribution of sub-focal times}

To complete the proof of Proposition~\ref{prop1.2}, we must remove the assumption that
at our self-focal point we have $\Phi_t(x,\xi)\notin S^*_xM$ for any $\xi$ if $t$
is not an integer multiple of the return time $\ell$.  As before, we may assume
that $\ell=1$.

We note that if $\delta>0$ is smaller than the injectivity radius of $(M,g)$, it
follows that we can never have $\Phi_t(x,\xi)\in S^*_xM$ for some
$t\in [\nu-\delta,\nu+\delta]$ with $\nu\in \Z$ since $\Phi_\nu(x,\xi)\in S^*_xM$.
Thus, in any such time interval, $(t,y)\to \bigl(e^{itP}\bigr)(y,y)$ is smooth
at $y=x$ if $t\in [\nu-\delta,\nu+\delta]\backslash \{\nu\}$.

To use this, fix $\beta\in C^\infty_0(\R)$ satisfying
\begin{equation}\label{z.1}
\beta(s)=1, \, \, \, |s|\le \frac34 \delta, \quad \text{and } \, \, \, 
\text{supp }\beta \subset (-\delta,\delta).
\end{equation}
We then can write
$$\rho\bigl(T(\la-P)\bigr)(y,y)=K(T,\la;y)+R(T,\la;y),$$
where
\begin{equation}\label{z.2}
K(T,\la;y)=\frac1{2\pi T}\sum_{|\nu|\le T}\int 
\beta(t-\nu) \, \Hat \rho(t/T) \, 
\bigl(e^{itP}\bigr)(y,y) \, e^{-it\la} \, dt,
\end{equation}
and
\begin{equation}\label{z.3'}
R(T,\la;y)=\frac1{2\pi T}\sum_{|\nu|\le T}\int 
\bigl(1-\beta(t-\nu)\bigr) \, \Hat \rho(t/T) \, 
\bigl(e^{itP}\bigr)(y,y) \, e^{-it\la} \, dt.
\end{equation}

Since for every $\nu\in \Z$, $|\nu|\le T$, 
$\beta(t-\nu) \bigl(e^{itP}\bigr)(y,y)$
can be written as in the right side of \eqref{parametrix} where the amplitude
satisfies \eqref{parametrix2} and \eqref{b.8}, it is clear that the proof of \eqref{1.11}
under the assumption that there are no subfocal times shows that, given $T=T(\e)$ large
enough and fixed, there is a neighborhood ${\mathcal N}_0$ of $x$ so that
\begin{equation}\label{z.3}
|K(T,\la; y)|\le \e \la^{n-1}+O_T(\la^{n-2}), \quad y \in {\mathcal N}_0.
\end{equation}
As a result, the proof of Proposition~\ref{prop1.2} would be complete if we could show
that 
for a given fixed $T$ there is a
a neighborhood ${\mathcal N}_1={\mathcal N}_1(\e,T)$ of $x$
so that
\begin{equation}\label{z.4}
|R(T,\la;y)|\le \e\la^{n-1}+O_T(\la^{n-2}),  \quad y\in {\mathcal N}_1.
\end{equation}
For if
we take ${\mathcal N}={\mathcal N}_0\cap {\mathcal N}_1$, we then have
$$\rho\bigl(T(\la-P)\bigr)(y,y)\le 2\e \la^{n-1} +O_T(\la^{n-2}), \quad 
y\in {\mathcal N},$$
which implies \eqref{1.11}, as $\e>0$ is arbitrary.

To prove \eqref{z.4}, we note that the set
$${\mathcal E}_x=\bigl\{(x,\xi)\in S^*_xM: \, \, 
\Phi_t(x,\xi)=(x,\eta), \, \, \text{some } \, \eta, \, \, 
t\in [-T,T] \,  \backslash \bigcup_{|\nu|\le T}(\nu-\delta/2, \nu+\delta/2)\bigr\}
$$
is closed.  Moreover, since we are assuming that $(M,g)$ is real analytic and we are avoiding the focal times
$\nu\in \Z$, this subset of $S^*_xM$ must have measure zero (see \cite[p. 416]{SZ}).  Therefore, working in local coordinates
in the cotangent bundle, given $\e_0>0$, we can choose a $C^\infty(\Rn \backslash 0)$ function $b(\xi)$ which is homogeneous of degree zero
so that if
$B(\xi)=1-b(\xi)$,
then
\begin{equation}\label{z.5}
0\le b\le 1, \quad \int_{S^{n-1}}b \, d\omega <\e_0, \quad {\mathcal E}_x\cap \supp B=\emptyset.
\end{equation}
The last condition means that if we write $\Phi_t(x,\xi)=(x(t),\xi(t))$ then if $\xi \in \supp B$
$$\gamma_\xi =\bigl\{x(t): \, \, t\in [-T,T] \,  \backslash \bigcup_{|\nu|\le T}(\nu-\delta/2, \nu+\delta/2)\bigr\}$$
is a union of geodesic segments all of which are disjoint from $x$.  Since the geodesic distance from $\gamma_\xi$ to $x$ is a continuous
function of $\xi\in S^*_xM$, it follows that
$$\min_{\xi\in \supp B}d_g(x,\gamma_\xi)>0.$$
Similarly, there must be a neighborhood ${\mathcal N}$ of $x$ such that for $y\in {\mathcal N}$ if
$(y(t),\xi(t))=\Phi_t(y,\xi)$ then
\begin{equation}\label{z.6}
\bigl\{y(t): \, \, t\in [-T,T] \,  \backslash \bigcup_{|\nu|\le T}(\nu-\delta/2, \nu+\delta/2)\bigr\}
\notin {\mathcal N}, \quad \text{if } \, \, \xi\in \supp B.
\end{equation}
Therefore, by H\"ormander's propagation of singularities theorem \cite{Hsing}, if
$\Psi\in C^\infty_0$ equals one near $x$ but is supported in ${\mathcal N}$, it follows that
if we let
$B(y,D)$ be the operator
with symbol 
$B(y,\xi)=\Psi(y)B(\xi)$,
then
$$\sum_{|\nu|\le T}\bigl(1-\beta(t-\nu)) \, \bigl(B(y,D)\circ e^{itP}\bigr)(y,y)\in C^\infty([-T,T]\times M),$$
due to the fact that $\beta(t-\nu)$ equals one on $[\nu-\frac34\delta,\nu+\frac34\delta]$.
Consequently,
$$\frac1{2\pi T}\sum_{|\nu|\le T}\int 
\bigl(1-\beta(t-\nu)\bigr) \, \Hat \rho(t/T) \, 
\bigl(B(y,D)\circ e^{itP}\bigr)(y,y) \, e^{-it\la} \, dt=O_{T,B}(1).$$

Since, if $b(y,D)$ is the operator with symbol $\Psi(y)b(\xi)$, we have
\begin{multline*}\Psi(y)R(T,\la;y)
\\
=\frac1{2\pi T}\sum_{|\nu|\le T}\int 
\bigl(1-\beta(t-\nu)\bigr) \, \Hat \rho(t/T) \, 
\Bigl((B(y,D)+b(y,D)) \, \circ \, e^{itP}\Bigr)(y,y) \, e^{-it\la} \, dt,
\end{multline*}
we would therefore obtain \eqref{z.4} if we could show that
\begin{multline}\label{z.7}
\frac1{2\pi T}\Bigl|\int 
\sum_{|\nu|\le T}\bigl(1-\beta(t-\nu)\bigr) \, \Hat \rho(t/T) \, 
\bigl(b\circ e^{itP}\bigr)(y,y) \, e^{-it\la} \, dt\Bigr|
\\
\le \e\la^{n-1}+O_{T,b}(\la^{n-2}).
\end{multline}

If 
$$m_{T,\beta}(\tau)=\sum_{|\nu|\le T}
\int \bigl(1-\beta(t-\nu)\bigr)\Hat \rho(t/T) e^{-it\tau}\, d\tau,$$
then the quantity in the left side of \eqref{z.7} is
$$\left| \sum_{j=0}^\infty m_{T,\beta}(\la-\la_j) \, (be_j)(y) \, e_j(y)\right|.$$
Since $m_{T,\beta}\in {\mathcal S}(\R)$, it follows that
$$|m_{T,\beta}(\tau)|\le C_{T,\beta,N}(1+|\tau|)^{-N}, \quad N=1,2,3,\dots .$$
As a result, the left side of \eqref{z.7} is bounded by
\begin{multline*}C_{T,\beta,N}\sum_{j=0}^\infty (1+|\la-\la_j|)^{-N}|be_j(y)| \, |e_j(y)|
\\
\le C_{T,\beta,N}\Bigl(\sum_{j=0}^\infty (1+|\la-\la_j|)^{-N}|be_j(y)|^2\Bigr)^{\frac12} \, \Bigl(\sum_{j=0}^\infty (1+|\la-\la_j|)^{-N}|e_j(y)|^2\Bigr)^{\frac12},
\end{multline*}
using the Cauchy-Schwarz inequality in the last step.
Recall that, by the generalized local Weyl formula (see e.g. \cite[Theorem 5.2.3]{SHang}), there is a constant $C$ depending only on $(M,g)$ so that if $A$ is a 
classical zero order pseudodifferential operator with principal symbol $a(x,\xi)$,
we have
\begin{equation}\label{z.8}
\sum_{\la_j\in [\la,\la+1]}|Ae_j(y)|^2
\le C\la^{n-1}\int_{p(y,\xi)\le 1}|a(y,\xi)|^2 \, d\xi +O_A(\la^{n-2}).
\end{equation}
Combining this with the preceding inequality implies that the left side of
\eqref{z.7} is
$$\le C_{M,T,\beta}\la^{n-1}\Bigl(\int_{S^{n-1}}|b|^2 \, d\omega\Bigr)^{\frac12}
+O_{M,T,\beta,b}(\la^{n-2}),$$
and, therefore, if $\e_0$ in \eqref{z.5} is small enough, we obtain \eqref{z.7},
which completes the proof.

\subsection{Analysis near non-focal points}

We shall now give the proof of Proposition~\ref{prop1.3}.  It is very similar to the argument that we just gave to prove bounds for the contribution of subfocal times for the estimate near self-focal points.

We are assuming that at a given $x\in M$ we now have $|{\mathcal L}_x|=0$, and
we need to show that \eqref{1.11} is valid.  Recall that ${\mathcal L}_x\subset S^*_xM$
is the set of unit directions such that $\Phi_t(x,\xi)\in S^*_xM$ for some $t\ne 0$.
In other words, the set of initial unit directions for geodesic loops through $x$.

Note that we then have that
$${\mathcal L}^T_x=\bigl\{\xi\in S^*_xM: \, \, 
\Phi_t(x,\xi)\in S^*_xM \, \, \text{for some } \, 
t\in [-T,T]\backslash \{0\}\bigr\}$$
is closed and of measure zero.  It follows that, given $\e_0>0$, we can find a
$b\in C^\infty(\Rn\backslash 0)$ which is homogeneous of degree zero so that we have the following analog of \eqref{z.5}
\begin{equation}\label{3.1}
0\le b\le 1, \quad \int_{S^{n-1}}b\, d\omega <\e_0, \quad
{\mathcal L}_x^T \cap \supp \beta = \emptyset.
\end{equation}
The last condition means that if $\xi\in \supp B$ and $\delta>0$ the geodesic two segments given by
$$\gamma_\xi = \{x(t): \, \, \delta\le |t|\le T\}$$
are disjoint from $x$ if $\Phi_t(x,\xi)=(x(t),\xi(t))$.  As before, this yields
$$\min_{\xi\in \supp B}d_g(x,\gamma_\xi)>0,$$
and so there must be a neighborhood ${\mathcal N}$ of $x$ so that if $y\in {\mathcal N}$ and $(y(t),\xi(t))=\Phi_t(y,\xi)$ then
$$\{y(t): \, \delta\le |t|\le T\}\notin {\mathcal N}, \quad \text{if } \, \, 
\xi\in \supp B.$$
Therefore, if $\Psi\in C^\infty_0$ equals one near $x$ but is supported in
${\mathcal N}$ and if we let $B(y,D)$ be the operator with symbol
$B(y,\xi)=\Psi(y)B(\xi)$, by propagation of singularities,
\begin{equation}\label{3.2}
\bigl(B\circ e^{itP}\bigr)(y,y)\in C^\infty(\{\delta \le |t|\le T\}\times M).
\end{equation}

To use this, fix $\beta\in C^\infty_0(\R)$ which equals one on $[-2\delta,2\delta]$.  Since for
$T\ge 1$ we have the uniform bounds
$$\left|\frac1{2\pi T}\int \beta(t)\rho(t/T) \, e^{it\tau}\, dt\right|
\le C_NT^{-1}(1+|\tau|)^{-N}, \quad N=1,2,3,\dots,$$
it follows from the sharp local Weyl law that
\begin{multline*}\left|\frac1{2\pi T}\int \beta(t)\rho(t/T) 
\bigl(e^{itP}\bigr)(y,y) \, e^{-it\la}\, dt\right|
\\
\le CT^{-1}\sum_{j=0}^\infty (1+|\la-\la_j|)^{-n-1}(e_j(y))^2
\le CT^{-1}\la^{n-1}.\end{multline*}
Therefore, by \eqref{rho}, and the fact that if $b(y,\xi)=\Psi(y)b(\xi)$, then
$B(y,D)+b(y,D)=\Psi(y)$, with $\Psi$ equal to one near $x$, we conclude that we
would have \eqref{1.11} if we could show that for fixed $T\ge1$ we have
\begin{equation}\label{3.3}
\frac1{2\pi T}
\int \bigl(1-\beta(t)\bigr) \rho(t/T) \, \bigl(B\circ e^{itP}\bigr)(y,y) \, 
e^{-it\la} \, dt =O_{B,T}(1),
\end{equation}
and
\begin{equation}\label{3.4}
\left|\frac1{2\pi T}\int \bigl(1-\beta(t)\bigr) \rho(t/T) \, \bigl(b\circ e^{itP}\bigr)(y,y) \, 
e^{-it\la} \, dt\right|\le C_T\sqrt{\e_0}\, \la^{n-1}+O_{b,T}(\la^{n-2}).
\end{equation}

The first bound \eqref{3.3}, just follows from \eqref{3.2} and the fact that
$\beta$ equals one on $[-2\delta,2\delta]$.

If 
$$m_{T,\beta}(\tau)=\frac1{2\pi T} \int
\bigl(1-\beta(t)\bigr) \, \rho(t/T) \, e^{it\tau} \, dt,$$
then the left side of \eqref{3.4} equals
$$\Bigl|\sum_{j=0}^\infty m_{T,\beta}(\la-\la_j) \, (be_j)(y) \, e_j(y)\Bigr|.$$
Since $m_{T,\beta}\in {\mathcal S}(\R)$, we can use the Cauchy-Schwarz inequality
to see that for every $N=1,2,3,\dots$ this is dominated by a constant
depending on $T$, $\beta$ and $N$ times
$$\bigl(\sum_{j=0}^\infty (1+|\la-\la_j|)^{-N}|be_j(y)|^2\bigr)^{\frac12} \, 
\bigl(\sum_{j=0}^\infty (1+|\la-\la_j|)^{-N}|e_j(y)|^2\bigr)^{\frac12}.$$
Thus, \eqref{3.4} follows from \eqref{z.8}, which completes the proof of
Proposition~\ref{prop1.3}.

\section{$\Omega(\la^{\frac{n-1}2})$ bounds at self-focal points with 
$U$-invariant functions}

In this section we shall finish the proof of Theorem~\ref{theorem1.1}.  Since we have just shown that we have \eqref{1.o} when at every self-focal point $x\in M$ the associated Perron-Frobenius operator $U=U_x$ in \eqref{1.6} has no nonzero $L^2(S^*_xM)$-invariant
functions, we would be done if we could establish the following.

\begin{proposition}\label{propm}
Let $(M,g)$ be a compact boundaryless real analytic Riemannian manifold of dimension
$n\ge 2$.  Assume that $x$ is a self-focal point and that $\ell>0$ is the first
return time for the geodesic flow.  Suppose further that, if $U: L^2(S^*_xM)\to 
L^2(S^*_xM)$ is the associated Perron-Frobenius map over $x$, we have
\begin{equation}\label{m.1}
Ug=g, \quad \text{some } \, \, 0\ne g\in L^2(S^*_xM).
\end{equation}
Let $\beta$ denote the number of conjugate points counted with multiplicity
along $\gamma(t)$, $0<t\le \ell$, with $\gamma(t)$ being a unit speed geodesic
starting at $x$.  Then if
\begin{equation}\label{m.2}
\mu_k=\frac{2\pi}\ell\bigl(k+\beta/4\bigr),
\end{equation}
there is a $c=c(M)>0$ so that
\begin{equation}\label{m.3}
\liminf_{k\to \infty}\mu_k^{-(n-1)}\Bigl(\, \sum_{\la_j\in [\mu_k,\mu_k+\delta]}
\bigl((e_j(x)\bigr)^2\, \Bigr)\ge c, \quad \text{if } \, \delta>0.
\end{equation}
\end{proposition}

The number $\beta$ here is independent of the geodesic starting at $x$, and it is
also commonly referred to as the Maslov index of the geodesic.

Since
$$\bigl\|\chi_{[\mu_k,\mu_k+\delta]}\bigr\|^2_{L^2(M)\to L^\infty(M)}=
\sup_{y\in M} \sum_{\la_j\in [\mu_k,\mu_k+\delta]} \bigl(e_j(y)\bigr)^2,
$$
it is clear that if \eqref{m.3} is valid then we cannot have \eqref{1.o2}, which is
the remaining part of Theorem~\ref{theorem1.1}.

\medskip

To prove the Proposition, as above, let $\rho\in {\mathcal S}(\R)$ be as
in \eqref{1.8}.  Then clearly, we would have \eqref{m.3} if we could show 
that there is a {\em uniform} constant $c>0$ so that whenever $T\gg 1$ is fixed we have
\begin{equation}\label{m.4}
\sum_{j=0}^\infty \rho\bigl(T(\mu_k-\la_j)\bigr) \, 
\bigl((e_j(x)\bigr)^2 \ge c \mu_k^{n-1}, \quad k\ge N_T,
\end{equation}
for some $N_T<\infty$.
One obtains \eqref{m.3} from this by taking $T=\delta^{-1}$ after recalling
that $\rho\ge 0$ and $\rho(0)=1$.

To prove \eqref{m.4}, as before, we may assume that $\ell=1$.  For the sake of
simplicity, we shall also assume that there are no sub-focal times, since
by the argument at the end of \S 2.1, their contributions to 
\eqref{m.4} will be $o(\nu_k^{n-1})$.  We therefore, are assuming
that 
$$\Phi_t(\xi, \xi)\ne (x,\xi), \, \, \forall \xi\in S^*_xM, 
\quad \text{if } \, \, t \notin {\mathbb Z}.
$$
Then, by \eqref{b.18}, we have that
\begin{multline*}
(2\pi)^{n}\sum_{j=0}^\infty \rho\bigl(T(\mu_k-\la_j)\bigr) \, 
\bigl((e_j(x)\bigr)^2
\\
= \mu_k^{n-1} T^{-1}\sum_{\nu=-\infty}^\infty\int  e^{-i\nu \mu_k}
\hat \rho(\nu/T) \sum_{j=1}^{N(\nu)}a_{\nu,j}(\nu,x,x,\mu_k\omega)\, d\omega
+O_T(\mu_k^{n-2}).
\end{multline*}
Since (see \S 3 of \cite{DG}, (3.2.15) in \cite{HFIO} and \cite{Be}) we have the more
precise version of \eqref{b.20}
\begin{equation}\label{m.5}\sum_{j=1}^{N(\nu)}a_{\nu,j}(\nu,x,x,\mu_k \omega)=i^{\nu\beta} \sqrt{J^\nu(\omega)} +O_T(\mu_k^{-1}),
\end{equation}
if $\Hat \rho(\nu/T)\ne 0$, and since $\mu_k$ is given by \eqref{m.2} with
$\ell =1$, we conclude that, modulo $O_T(\mu_k^{n-2})$ terms, the left
side of \eqref{m.4} equals $(2\pi)^{-n}$ times
$$\mu_k^{n-1} \Bigl(\, T^{-1}\sum_{\nu=-\infty}^\infty 
\int \hat \rho(\nu/T) \sqrt{J^\nu(\omega)} \, d\omega \, \Bigr)
=\mu_k^{n-1} \Bigl( \, T^{-1} \sum_{\nu=-\infty}^\infty \hat \rho(\nu/T) \, \langle U^\nu \1, \1\rangle
\, \Bigr).$$
Since $\hat \rho$ is nonnegative and $\hat \rho(0)>0$, it follows that there
is a constant $c_0>0$, which is independent of $T$, so that
$$T^{-1} \sum_{\nu=-\infty}^\infty \hat \rho(\nu/T) \, \langle U^\nu \1, \1\rangle
\ge c_0 \, M^{-1}\sum_{\nu=-M}^M\langle U^\nu \1, \1\rangle, \quad M=c_0T.$$
By von Neumann's ergodic theorem we have that
$$\lim_{M\to \infty}\,\frac1{2M}\sum_{\nu=-M}^M\langle U^\nu \1, \1\rangle
=\langle \varPi(\1), \1\rangle,$$
where $\varPi(1)$ denotes the projection of $\1$ onto the $U$-invariant
subspace of $L^2(S^*_xM)$.  Since our assumption \eqref{m.1} gives that
$$\langle \varPi(\1), \1\rangle>0,$$
we conclude that \eqref{m.4} must be valid, which completes the proof
of Proposition~\ref{propm}.

\section{Bounds for quasi-modes}

In this section, we shall prove Corollary~\ref{corr}.  Recall that we are assuming that
$\phi_\la$ are quasi-modes satisfying
\eqref{q.2}, i.e.,
\begin{equation}\label{4.2}
\int |\phi_\la|^2 \, dV=1, \quad
\text{and } \, \, \|S_{[2\la,\infty)}\phi_\la \|_{L^\infty(M)}+
\|(\Delta+\la^2)\phi_\la\|_{L^2(M)}=o(\la),
\end{equation}
with $S_{[2\la,\infty)}$ denoting the projection on to frequencies in $[2\la,\infty)$.
By Theorem~\ref{theorem1.1}, if there is a 
a self-focal
point $x$ for which the operators in \eqref{1.6} have a nontrivial invariant function
satisfying \eqref{m.1}, then we know that there is a uniform constant $c>0$ so that, given any
$\delta>0$ we have
$$\limsup_{\la\to \infty}\la^{-\frac{n-1}2}\bigl\| \chi_{[\la,\la+\delta]}\bigr\|_{L^2(M)\to L^\infty(M)}\ge c.$$
Thus, in this case, if $\delta_j\to 0$, we can find $\phi_{\la_j}$ with spectrum in $[\la_j,\la_j+\delta]$ satisfying
$\|\phi_{\la_j}\|_2=1$ and $\|\phi_{\la_j}\|_\infty \ge c\la_j^{\frac{n-1}2}$.  Consequently since
the $\{\phi_{\la_j}\}$ are clearly quasi-modes of order $o(\la)$, to finish the proof we just need
to prove that \eqref{q.3} is valid when all the self-focal points are dissipative.

To do this, we just use the fact that,
by Lemma 2.5 in \cite{STZ}, we must have
$$\|\phi_\la\|_{L^\infty(M)}=o(\la^{\frac{n-1}2})$$
in this case, which completes the proof.

\begin{rem}  
It would be interesting to see whether we can obtain $\Omega(\mu_k^{\frac{n-1}2})$ bounds for quasi-modes of {\em order zero} whenever
there is a self-focal point $x\in M$ satisfying \eqref{m.1}.  We recall these are ones satisfying
the stronger variant of \eqref{q.2}, which says that
$$\|(\Delta+\mu_k^2)\phi_{\mu_k}\|_{L^2(M)}=O(1), \quad \text{and } \, \, \int |\phi_{\mu_k}|^2 \, dV=1.$$
Such quasi-modes were constructed in \cite{STZ} when the map $\xi\to \eta(\xi)$ equals the identity map on an open subset of
$S^*_xM$, which, of course, is a stronger condition than \eqref{m.1}. It is plausible that in the real analytic setting,
one may construct quasi-modes of order zero corresponding to a self-focal point such that $U_x$ has an invariant
$L^2$ function. The idea is to use $|f|^2 d\mu_x$ as the `symbol' or invariant measure in the quasi-mode construction.
In the real analytic setting, one can try to `quantize' $(\Lambda_x, |f|^2 d\mu_x)$ using Toeplitz quantization 
in a Grauert tube around $M$; we plan to investigate this in future work. 
\end{rem}

\end{document}